\documentclass[11pt]{article}
\usepackage{amssymb}
\def\bbb#1{\hbox {{\gordas #1}}}

\font\gordas = msbm10 at 11pt
\def\bbb#1{\hbox {{\gordotas #1}}}
\def\erre{{\bbb R}}

\font\gordas = msbm10 at 12pt
\def\bbb#1{\hbox {{\gordas #1}}}

\def\UNO{1\mkern-7mu1}

\newtheorem{theorem}{Theorem}[section]

\newtheorem{proposition}[theorem]{Proposition}
\newtheorem{definition}[theorem]{Definition}
\newtheorem{remark}[theorem]{Remark}

\begin{document}
\noindent
%extfBm (16/04/07)
\begin{center}
{\bf\Large Some Extensions of Fractional Brownian Motion and \\[.25cm] Sub-Fractional Brownian Motion Related to Particle Systems}  \\[1cm]
T. BOJDECKI$^1$\\
{\it Institute of Mathematics, University of Warsaw\\
ul. Banacha 2, 02-097 Warsaw, Poland (e.mail: tobojd@mimuw.edu.pl)}\\[1cm]
L.G. GOROSTIZA$^2$\\
{\it Centro de Investigaci\'on y de Estudios Avanzados\\
07000 M\'exico, D.F., Mexico (e.mail: lgorosti@math.cinvestav.mx)}\\[1cm]
A. TALARCZYK$^1$\\
{\it Institute of Mathematics, University of Warsaw\\
ul. Banacha 2, 02-097 Warsaw, Poland (e.mail: annatal@mimuw.edu.pl)}\\[1cm]
\end{center}
\footnote{\kern -.6cm $^1$ Research supported by MNiSW grant 1P03A01129 (Poland).\\
$^2$ Research supported by CONACyT grant 45684-F (Mexico).}
\begin{abstract}

\vglue .25cm
In this paper we study three self-similar, long-range dependence, Gaussian processes. The first one, with covariance
$$
\int^{s\wedge t}_0 u^a [(t-u)^b+(s-u)^b]du,
$$
parameters $a>-1$, $-1<b\leq 1$, $|b|\leq 1+a$, corresponds to fractional Brownian motion for $a=0$, $-1<b<1$. The second one, with covariance
$$
(2-h)\biggl(s^h+t^h-\frac{1}{2}[(s+t)^h +|s-t|^h]\biggr),
$$
parameter $0<h\leq 4$,  corresponds to sub-fractional Brownian motion for $0<h<2$. The third one, with covariance 
$$
-\left(s^2\log s + t^2\log t -\frac{1}{2}[(s+t)^2 \log (s+t) +(s-t)^2 \log |s-t|]\right),
$$ 
is related to the second one. These processes come from occupation time fluctuations of certain particle systems for some values of the parameters.
\end{abstract}
\vglue .25cm
\noindent
\newpage
\section{Introduction}
\setcounter{section}{1}
\setcounter{equation}{0}
\label{sec:1}
Some types of self-similar, long-range dependence, centered Gaussian processes have been obtained from occupation time fluctuation limits of certain particle systems (they represent the temporal structures of the limits). These processes have parameters whose ranges of possible values are determined by properties of the particle systems. Since they are interesting processes in themselves, it is natural to ask if the ranges of values of the parameters can be extended and what the maximal ranges are, what properties the processes possess, and if the extended processes may also be related to some  stochastic  models. In this paper we will consider the first two questions for some extensions of fractional Brownian motion (fBm) and sub-fractional Brownian motion (sfBm). The third question, which is of a different kind (weak convergence of processes related to the particle systems), is discussed in \cite{BGT6}.

In the remainder of the Introduction we will discuss the motivations behind the present article in connection with particle systems. In the following sections we will study the processes; most of the contents in those sections  can be read independently of this  discussion.

We will study a centered Gaussian process $(\xi_t)_{t\geq 0}$ with covariance $Q(s,t)=E\xi_s\xi_t$, $s,t\geq 0$, given by
\begin{equation}
\label{eq:1.1}
Q(s,t)=\int^{s\wedge t}_0 u^a [(t-u)^b+(s-u)^b]du.\end{equation}
For $a=0$, (\ref{eq:1.1}) becomes 
\begin{equation}
\label{eq:1.2}
\frac{1}{b+1}(t^{b+1}+s^{b+1}-|s-t|^{b+1}),
\end{equation}
which for $-1<b<1$  corresponds to the covariance of the well-known fBm with Hurst parameter $\frac{1}{2}(b+1)$, Brownian motion (Bm) for $b=0$. We call the process $\xi$ with covariance (\ref{eq:1.1}) a {\it weighted fractional Brownian motion} (wfBm), refering to the weight function $u^a$ in the integral. Thus, wfBm's are a family of processes which extend fBm's.

The process $\xi$ appeared in \cite{BGT5} in a limit of occupation time fluctuations of a system of independent particles moving in $\erre^d$ according a symmetric $\alpha$-stable L\'evy process, $0<\alpha \leq 2$, started from an inhomogeneous Poisson configuration with intensity measure $dx/(1+|x|^\gamma)$, $0<\gamma \leq d=1<\alpha$, $a=-\gamma /\alpha$, $b=1-1/\alpha$, the ranges of values of $a$ and $b$ being $-1<a<0$ and $0<b\leq 1+a$. The homogeneous case $(\gamma=0$, fBm) was studied in \cite{BGT2, BGT3}. The process $\xi$ also appears  in \cite{BGT6} in a high-density limit of occupation time fluctuations of the above mentioned particle system, where the initial Poisson configuration has finite intensity measure, with $d=1<\alpha$, $a=-1/\alpha$, $b=1-1/\alpha$. From these results it follows that for $-1<a<0$ and any $0<b\leq (1+a)\wedge (1/2)$, the  kernel (\ref{eq:1.1}) is positive-definite, because it is a limit of covariances. The question that interests us here is what are the maximal ranges of values of the parameters $a$ and $b$ for which (\ref{eq:1.1}) is a positive-definite kernel (the conditions $a>-1$ and $b>-1$ are obviously necessary). We will show that the maximal ranges are $a>-1$,  $-1<b\leq 1$, and $|b|\leq 1+a$, and we will study some properties of the Gaussian process $\xi$ with covariance (\ref{eq:1.1}). These properties, which are analogous to those of fBm, are self-similarity, path continuity, behavior of increments,  long-range dependence, non-semimartingale, and others. Since wfBm widens the scope of behaviour of fBm, perhaps it may be useful in some applications. 

A different extension of fBm is the {\it bi-fractional Brownian motion} (bfBm) introduced in \cite{HV}, with covariance
\begin{equation}
\label{eq:1.3}
(s^{2H}+t^{2H})^k -|s-t|^{2Hk},
\end{equation}
where $0<H<1$ and $0<k\leq 1$, corresponding to fBm with Hurst parameter $H$ for $k=1$. This process was obtained analytically from the covariance of fBm, not from a specific stochastic model, whereas  the process with covariance (\ref{eq:1.1}) appeared first in connection with a particle system. The process  bfBm is also studied in \cite{RT,TX}. With $H=1/2$, bfBm can be extended for $1<k<2$, and in this way it  plays a role in our next process, which we describe presently.

Let $(\zeta_t)_{t\geq 0}$ be a centered Gaussian process with covariance $K(s,t)=E\zeta_s \zeta_t$, $s,t\geq 0$, given by 
\begin{equation}
\label{eq:1.4}
K(s,t)=(2-h) \biggl(s^h +t^h-\frac{1}{2}[(s+t)^h +|s-t|^h]\biggr),
\end{equation}
where $0<h\leq 4$, which for $h=1$ corresponds to Bm. This process was found in \cite{BGT2, BGT3} for $1<h\leq 2$ (where  the covariance was written without the factor $(2-h)$, which is irrelevant there) in a limit of occupation time fluctuations of a particle system, which is as the one described above in the homogeneous case $(\gamma =0)$, and in addition the particles undergo critical binary branching, i.e., each particle independently, at an exponentially distributed lifetime, disappears with probability $1/2$ or is replaced by two particles at the same site with probability $1/2$. (See also \cite{I}   in the setting or superprocesses with $\alpha =2$.) In this model $d$ is restricted to $\alpha <d<2\alpha$, and $h=3-d/\alpha$; hence the range of values of $h$ for the particle model is $1<h<2$. However, existence of the process $\zeta$ for all $0<h<2$ follows from the equality in distribution
\begin{equation}
\label{eq:1.5}
(\zeta_t)_{t\geq 0}=c(\xi_t+\xi_{-t})_{t\geq 0},
\end{equation}
where $\xi$ is fBm defined for all $-\infty <t<\infty$ and $c=1/\sqrt{2}$  \cite{BGT2}. The process $\zeta$ with $0<h<2$ was called {\it sub-fractional Brownian motion} in \cite{BGT2} because it has properties analogous to those of fBm, and it is  intermediate between Bm and fBm in the sense that its increments on non-overlapping intervals are more weakly correlated and their covariance decays faster than for fBm. Properties of sfBm have been studied in \cite{BGT2} and \cite{Tu1,Tu2}. The process (\ref{eq:1.5}) with $c=1/2$ is called the even part of fBm in \cite{DZ}, where it appears in a different context.

We will show that $K(s,t)$ is a covariance also for $2<h\leq 4$ (degenerate for $k=4$). Since the coefficient $(2-h)$ in  $K(s,t)$ changes sign from positive to negative as $h$ passes from below to above $2$, we call the process $\zeta$ for $2<h<4$ a {\it negative sub-fractional Brownian motion} (nsfBm).

Note that in contrast to this extension of sfBm, in the case of fBm it is not possible to extend the covariance for Hurst parameter $>1$ (the kernel is not positive definite).

It turns out that nsfBm for any $2<h<4$ is represented as the integral of the extended bfBm referred to above for $k=h-2$  (with $H=1/2$), i.e., with covariance
\begin{equation}
\label{eq:1.6}
(s+t)^{h-2}-|s-t|^{h-2}.
\end{equation}
This implies that nsfBm is a semimartingale, as opposed to sfBm, which is not. We will also discuss other properties of nsfBm.

For $2<h<5/2$, the covariance (\ref{eq:1.4}) is a limit of covariances of occupation time fluctuations of the branching particle system described above in the homogeneous case for $d=1<\alpha$, with $h=3-d/\alpha$ (therefore it {\it is} a covariance), but this covariance does not correspond to the occupation time fluctuation limit of the particle system, because the system becomes extinct locally for $d\leq \alpha$ \cite{GW}. For $d=\alpha$ there is a functional ergodic theorem  \cite{T}. On the other hand,  the  high-density occupation time fluctuation limits of the branching particle system obtained in \cite{BGT6} yield both sfBm for $\alpha<d<2\alpha$ and nsfBm for $d<\alpha$. In this setup, the case $d=\alpha$ leads to a centered Gaussian process $(\eta_t)_{t\geq 0}$ with covariance $R(s,t)=E\eta_s\eta_t$, $s,t\geq 0$, given by
\begin{equation}
\label{eq:1.7}
R(s,t)=-\biggl(s^2\log s+t^2\log t-\frac{1}{2}[(s+t)^2 \log (s+t)+(s-t)^2\log |s-t|]\biggr)
\end{equation}
\cite{BGT6}, and we will also study this process. Although it has a similar interpretation to that of nsfBm regarding the particle system, it is not a semimartingale.

As we have seen, the processes described above arise in the limits of certain particle systems for some values of the parameters. It would be interesting to find  particle systems that produce these processes also for the remaining values of the parameters.

We remark that there exist more general models that yield, in the limit, stable, self-similar, long-range dependence processes, which in the Gaussian case coincide with the processes above for some of their parameter values \cite{BGT2, BGT3}. 

The present paper was written separately from \cite{BGT5, BGT6} because it can be read independently. The  results and the methods of proof in those papers refer to convergence of occupation times of particle systems.

Self-similar processes with long-range dependence have attracted much attention recently for their applications and their intrinsic mathematical interest (see e.g. \cite{DOT, Ta}). Therefore it seems worthwhile to investigate new processes of  this type, and stochastic models where they arise.

\section{Weighted fractional Brownian motion}
\label{sec:2}
\setcounter{equation}{0}
\subsection{Results}
\label{sub:2.1}
Let us recall the funtion $Q$ defined by (\ref{eq:1.1}):
\begin{equation}
\label{eq:2.1}
Q(s,t)=\int^{s\wedge t}_0 u^a[(t-u)^b +(s-u)^b]du, \quad s,t \geq 0.
\end{equation}

We start by discussing conditions under which this function (symmetric, continuous) is the covariance of a stochastic process.
\begin{theorem}
\label{T2.1} The function $Q$ defined by (\ref{eq:2.1}) is positive-definite if and only if $a$ and $b$ satisfy the  conditions
\begin{equation}
\label{eq:2.2}
a>-1, \quad -1<b\leq 1, \quad\quad |b|\leq 1+a.
\end{equation}
\end{theorem}

The proof of this theorem will be given in the next subsection.
\begin{definition}
{\rm The centered Gaussian process $\xi = (\xi_t)_{t\geq 0}$ with covariance function  (\ref{eq:2.1}) satisfying conditions (\ref{eq:2.2}) will be called} {\it weighted fractional Brownian motion (wfBm) with parameters $a$ and $b$.}
\end{definition}
\begin{remark}
\label{T2.3} {\rm (a) It is immediately seen that for $a=0$, $\xi$ is the usual fBm with Hurst parameter $\frac{1}{2}(b+1)$ (up to a multiplicative  constant).

\noindent
(b) Analogously to the case of fBm, in our definition of wfBm we exclude $b=1$. It is easy to see that for $b=1$ (and $a\geq 0$, according to Theorem \ref{T2.1}), (\ref{eq:2.1}) is the covariance  of the process
$$
\xi_t=\int^t_0 w_{r^a}dr,
$$
where $w$ is a standard Bm.

\noindent
(c) For $b=0$, $\xi$ is a time-inhomogeneous Bm.}
\end{remark}

In the next theorem we collect the main properties  of wfBm.

\begin{theorem}
\label{T2.4} The weighted fractional Brownian motion $\xi$ with parameters $a$ and $b$ has the following properties:

\noindent
(1) Self-similarity:
$$
(\xi_{ct})_{t\geq 0} \stackrel{d}{=} (c^{(1+a+b)/2}\xi_t)_{t\geq 0}\quad \hbox{\it for each} \quad c>0.
$$
(2) Second moments of increments: for $0\leq s<t$,
\begin{equation}
\label{eq:2.3}
E(\xi_t-\xi_s)^2=2 \int^t_s u^a (t-u)^b du,
\end{equation}
\begin{equation}
\label{eq:2.4}
E(\xi_t-\xi_s)^2 \leq C|t-s|^{b+1},
\end{equation}
if $a\geq 0$, $s,t\leq T$ for any $T>0$ with $C=C(T)$, and also if $a<0$, $s,t\geq \varepsilon >0$ for any $\varepsilon >0$, with $C=C(\varepsilon)$;
\begin{equation}
\label{eq:2.5}
E(\xi_t -\xi_s)^2 \leq C|t-s|^{1+a+b}, \quad s, t\geq 0,
\end{equation}
if $a<0$, $1+a+b>0$;
\begin{equation}
\label{eq:2.6}
E(\xi_t-\xi_s)^2 \geq C|t-s|^{b+1},
\end{equation}
if $a>0$, $s,t\geq \varepsilon >0$ for any $\varepsilon >0$ with $C=C(\varepsilon)$ and also if $a\leq 0$, $s, t\leq T$ for any $T>0$ with $C=C(T)$.

\noindent
(3) Path continuity: $\xi$ is a continuous process with the only exception of the case $a<0$, $b<0$, $a+b=-1$, where $\xi$ is discontinuous at $0$. 

\noindent
(4) Covariance of increments: For $0\leq r<v\leq s<t$,
\begin{equation}
\label{eq:2.7}
Q(r,v,s,t) = E((\xi_t-\xi_s) (\xi_v-\xi_r))
= \int^v_r u^a [(t-u)^b-(s-u)^b]du,
\end{equation}
hence 
$$
Q(r,v,s,t)\;\;\left\{
\begin{array}{lccr}
 & >0 & if & b>0,\\
 & =0 & if & b=0,\\
 & <0 & if & b<0.
\end{array}\right.
$$
(5) Long-range dependence: For $r,v,s,t$ as in (4),
\begin{equation}
\label{eq:2.8}
\lim_{T\rightarrow \infty} T^{1-b} Q(r,v,s+T, t+T) =\frac{b}{a+1} (t-s) (v^{a+1}-r^{a+1}).
\end{equation}
(6) Asymptotic homogeneity: The finite-dimensional distributions of the process $(T^{-a/2}(\xi_{t+T}-\xi_T))_{t\geq 0}$ converge as $T\rightarrow \infty$ to those of fBm with Hurst paramter $(1+b)/2$, multiplied by $(2/(1+b))^{1/2}$.

\noindent
(7) Short and long-time asymptotics:
\begin{eqnarray}
\label{eq:2.9}
\lim_{\varepsilon \searrow 0} \varepsilon^{-b-1} E(\xi_{t+\varepsilon}-\xi_t)^2 &=& \frac{2}{b+1} t^a,\\
\label{eq:2.10}
\lim_{T\rightarrow \infty} T^{-(1+a+b)} E(\xi_{t+T} -\xi_t)^2 &=& 2\int^1_0 u^a (1-u)^b du.
\end{eqnarray}
Hence $\xi$ has asymptotically stationary increments for long time intervals, but not for short time intervals.

\noindent
(8) $\xi$ is not a semimartingale if $b\neq 0$.

\noindent
(9) $\xi$ is not Markov if $b\neq 0$.
\end{theorem}

\subsection{Proofs}
\label{sub:2.2}
\noindent
{\bf Proof of Theorem \ref{T2.1}} Firstly, we prove positive definiteness of $Q$ in the case 
\begin{equation}
\label{eq:2.11}
a>-1, \quad -1<b \leq 0, \quad a+b+1\geq 0,
\end{equation}
(see (\ref{eq:2.2})). We have from (\ref{eq:2.1})
$$
Q(s,t)=Q_1(s,t) +Q_2 (s,t),
$$
where 
$$
Q_1(s,t) =\int^{s\wedge t}_0 u^a ((s\wedge t)-u)^b du, \quad Q_2 (s,t) =\int^{s\wedge t}_0 u^a ((s\vee t)-u)^bdu.
$$
It suffices to show that both $Q_1$ and $Q_2$  are positive-definite. 

$Q_1$ can be written as 
$$
Q_1(s,t) =(s\wedge t)^{1+a+b} \int^1_0 u^a (1-u)^b du,
$$
so it is positive definite by (\ref{eq:2.11}).

Next, since $b\leq 0$ we can write $Q_2$ as
$$Q_2 (s,t) =\int^\infty_0 u^a[(s-u)^b \wedge (t-u)^b]\UNO_{[0,s]}(u)\UNO_{[0,t]}(u)du,
$$
hence positive definiteness of $Q_2$ follows easily.

Now assume that
\begin{equation}
\label{eq:2.12}
a>-1, \quad 0<b\leq 1, \quad b\leq 1+a.
\end{equation}
From (\ref{eq:2.1}), $Q$ can be transformed into 

\begin{equation}
\label{eq:2.13}
Q(s,t)=b\int^s_0\int^t_0 (u\wedge r)^a |u-r|^{b-1}drdu,
\end{equation}
hence it is clear that $Q$ is positive-definite for $b=1$ (note that (\ref{eq:2.12}) implies $a\geq 0$ in this case).

Assume $b<1$; then, from (\ref{eq:2.13}),
$$
Q(s,t) =b \int^s_0 \int^t_0 (u\wedge r)^a (u\vee r)^{b-1}\biggl(1-\frac{u\wedge r}{u\vee r}\biggr)^{b-1}dudr.
$$
Using the power series expansion of $(1-\cdot)^{b-1}$  we obtain
$$
Q(s,t)=\sum^\infty_{n=0} b \frac{\Gamma (1-b +n)}{\Gamma (1-b)n!}\int^s_0 \int^t_0 (u\wedge r)^a (u\vee r)^{b-1}\frac{(u\wedge r)^n}{(u\vee r)^n}dudr.
$$
Each summand is positive definite, since
$$
(u\wedge r)^a (u\vee r)^{b-1}(u\wedge r)^n(u\vee r)^{-n}=u^ar^a (u^{b-1-a}\wedge r^{b-1-a})(u\wedge r)^n (u^{-n}\wedge r^{-n})
$$
(we have used that $b-1-a\leq 0$, by (\ref{eq:2.12})). Hence $Q$ is positive-definite in this case.

We will show that for the remaining values of parameters $Q$ is not positive-definite (recall that $Q$ is infinite if either $a$ or $b$ is $\leq -1$). More precisely, we will prove that there exists $t>0$ such that the covariance inequality
\begin{equation}
\label{eq:2.14}
Q(1,t)\leq (Q(1,1)Q(t,t))^{1/2}
\end{equation}
does not hold. Indeed, for $-1<b<0$, $a>-1$,  $a+b+1<0$ and $t\searrow 0$, the left-hand side of (\ref{eq:2.14}) is of order $t^{1+a+b}$ while the right-hand side is of order $t^{(1+a+b)/2}$.

For $a>-1$, $b>a+1$ and $t\nearrow \infty$ the left-hand side of (\ref{eq:2.14}) is of order $t^b$ and the right-hand side is of order $t^{(1+a+b)/2}$.

It remains to consider the case $1<b\leq a+1$. We show that (\ref{eq:2.14}) does not hold for $t=1+\varepsilon$, $\varepsilon \searrow 0$. Using convexity of the function $x^b$ for $b>1$, we have
\begin{eqnarray*}
Q(1, 1+\varepsilon) &\geq & 2 \int^1_0 u^a \biggl(1+\frac{\varepsilon}{2}-u\biggr)^b du\\
&=& 2\biggl(1+\frac{\varepsilon}{2}\biggr)^{a+b+1}
\int^{(1+\varepsilon/2)^{-1}}_0 u^a (1-u)^b du\\
&\geq & 2\biggl(1+\varepsilon +\frac{\varepsilon^2}{4}\biggr)^{(a+b+1)/2}
\biggl[\int^1_0 u^a(1-u)^b du -\biggl(\frac{\varepsilon}{2}\biggr)^{b+1}\frac{1}{b+1}\biggr].
\end{eqnarray*}
This implies
\begin{eqnarray}
\label{eq:2.15}
\lefteqn{Q(1,1+\varepsilon) -(Q(1,1)Q(1+\varepsilon, 1+\varepsilon))^{1/2}}\nonumber\\
&=&Q(1,1+\varepsilon)-2(1+\varepsilon)^{(a+b+1)/2}\int^1_0 u^a (1-u)^b du\nonumber\\
&\geq & 2\biggl[\biggl(\biggl(1+\varepsilon +\frac{\varepsilon^2}{4}\biggr)^{(a+b+1)/2}-(1+\varepsilon )^{(a+b+1)/2}\biggr)\int^1_0 u^a (1-u)^b du\biggr]\nonumber\\
&&- \frac{1}{2^b(b+1)}\biggl(1+\varepsilon +\frac{\varepsilon^2}{4}\biggr)^{(a+b+1)/2}\varepsilon^{b+1}
\geq  A\varepsilon^2 -B\varepsilon^{b+1}
\end{eqnarray}
for some positive constants $A,B$. In the last estimate we have used the fact that $a+b+1\geq 2b>2$.

The right-hand side of (\ref{eq:2.15}) is strictly positive for $\varepsilon$ sufficiently small, since $b>1$, so (\ref{eq:2.14}) does not hold for such $\varepsilon$.
\hfill $\Box$

\vglue .25cm
\noindent
{\bf Proof of Theorem 2.4}
\vglue .25cm
\noindent
(1) Self similarity follows immediately from (\ref{eq:2.1}).

\noindent
(2) Formula (\ref{eq:2.3}) is a direct consequence of (\ref{eq:2.1}), and (\ref{eq:2.4}) follows from (\ref{eq:2.3}).

To prove (\ref{eq:2.5}), first observe that if $a<0$, $b\geq 0$, (\ref{eq:2.3}) implies
$$
E(\xi_t-\xi_s)^2 \leq \frac{2}{a+1}|t-s|^b |t^{a+1}-s^{a+1}|\leq \frac{2}{a+1}|t-s|^{1+a+b}.
$$
Next assume that $a<0$, $b<0$, $1+a+b>0$. Fix any $\theta$ such that 
$|a|/(1+b)<\theta <1$, and put  $p=\theta/|a|$ and $q =\theta/(\theta -|a|)$. For $s<t$, the H\"older inequality applied to (\ref{eq:2.3}) yields
\begin{eqnarray*}
E(\xi_t-\xi_s)^2 &\leq & 2 \biggl(\int^t_s u^{ap}du\biggr)^{1/p}\biggl(\int^t_s(t-u)^{bq}du\biggr)^{1/q}\\
&\leq & C(t^{ap+1}-s^{ap+1})^{1/p}(t-s)^{(bq+1)/q}\\
&\leq & C(t-s)^{1+a+b},
\end{eqnarray*}
since $0 <ap+1<1$.

The inequality (\ref{eq:2.6}) follows from (\ref{eq:2.3}).

\noindent
(3) Path continuity of $\xi$ is a consequence of (\ref{eq:2.4}) and (\ref{eq:2.5}). If $a<0$, $b<0$, $a+b=-1$, then
$$
E\xi^2_t =2\int^1_0 u^a (1-u)^bds\quad {\rm for}\quad t>0.
$$
(4) Formula (\ref{eq:2.7}) follows from (\ref{eq:2.1}).

\noindent
(5) The limit (\ref{eq:2.8}) can be easily obtained from (\ref{eq:2.7}).

\noindent
(6) For $0\leq s\leq t$, by (\ref{eq:2.3}) and (\ref{eq:2.7}) we have 
\begin{eqnarray*}
\frac{1}{T^a} E((\xi_{t+T}-\xi_T) (\xi_{s+T} -\xi_T)) &=& \frac{1}{T^a} \int^{s+T}_T u^a ((t+T-u)^b +(s+T-u)^b)du\\
&=& \int^s_0 \biggl(\frac{u}{T}+1\biggr)^a ((t-u)^b +(s-u)^b)du\\
&\rightarrow & \frac{1}{b+1}(t^{b+1}+s^{b+1}-(t-s)^{b+1})\quad {\rm as}\quad T\rightarrow \infty,
\end{eqnarray*}
hence the assertion follows because $\xi$ is Gaussian.

\noindent
(7) The proofs of (\ref{eq:2.9}) and (\ref{eq:2.10}) are very easy, so we omit them.

\noindent
(8) The non-semimartingale property follows from (\ref{eq:2.4}), (\ref{eq:2.6}) and Corollary 2.1 in \cite{BGT1}.

\noindent
(9) It is not hard to see that for $b\neq 0$ the covariance (\ref{eq:2.1}) does not have the triangular property, so $\xi$ is not Markovian (see \cite{K}, Proposition 13.7). \hfill $\Box$

\section{Negative sub-fractional Brownian motion}
\setcounter{equation}{0}
\subsection{Results}

Let us recall the function $K$ defined by (\ref{eq:1.4}):
\begin{equation}
\label{eq:3.1}
K(s,t)=(2-h)\biggl(s^h+t^h-\frac{1}{2}[(s+t)^h+|s-t|^h]\biggr).
\end{equation}
We know that for $0<h<2$ this is the covariance of a sub-fractional Brownian motion \cite{BGT2}. Now we consider the case $h>2$.

\begin{theorem}
\label{T:3.1} 
For $2<h\leq 4$ there exists a centered Gaussian process $\zeta=(\zeta_t)_{t\geq 0}$ with covariance (\ref{eq:3.1}). This process has a representation
\begin{equation}
\label{eq:3.2}
\zeta_t =\biggl(\frac{1}{2}h(h-1)(h-2)\biggr)^{1/2}\int^t_0\vartheta_sds,
\end{equation}
\end{theorem}
{\it where $\vartheta$ is a centered continuous Gaussian process with covariance}
\begin{equation}
\label{eq:3.3}
K_0(s,t)=(s+t)^{h-2}-|s-t|^{h-2}.
\end{equation}
\begin{definition}
\label{D:3.2}
{\rm The process $\zeta$ described in Theorem 3.1 for $2<h<4$ will be called} {\it negative sub-fractional Brownian motion  (nsfBm) with parameter $h$}.
\end{definition}
\begin{remark}
\label{R:3.3}
{\rm (a) It is not difficult to see that (\ref{eq:3.1}) is not positive-definite if $h>4$. Indeed, for $t\to\infty,\quad K(1,t)$ is of order $t^{h-2}$, while $K(t,t)^{1/2}$ is of  order $t^{h/2}$.}

\noindent
{\rm (b) For $h=4,\zeta_t=ct^2\gamma$, where $\gamma$ is a standard Gaussian random variable and $c$ is a constant.}

\noindent
{\rm (c) It is easy to see by a simple covariance calculation that the process $\vartheta$ appearing in Theorem 3.1 is equal in distribution to the process $(\xi_t-\xi_{-t})_{t\geq 0}$, where $\xi$ is fBm with Hurst parameter $(h-2)/2$ (defined for all $-\infty <t<\infty$). The process $\frac{1}{2}(\xi_t-\xi_{-t})_{t\geq 0}$ is called the odd part of fBm in \cite{DZ}.}

\noindent
{\rm (d) For $2<h<4$ the covariance of nsfBm can be written as}
\begin{equation}
\label{eq:3.4}
K(s,t)=h(h-1)(h-2)^2\int^{s\wedge t}_0\int^s_r\int^t_r(u+u'-2r)^{h-3}du'dudr.
\end{equation}
\end{remark}

In the next proposition we gather basic properties of nsfBm.

\begin{proposition}
\label{P:3.4}
The negative sub-fractional Brownian motion $\zeta$ with parameter $h$ has the following properties:

\noindent
(1) $\zeta$ has continuously differentiable paths, hence it is a semimartingale.

\noindent
(2) Self-similarity:
$$(\zeta_{ct})_{t\geq 0}\stackrel{d}{=}(c^{h/2}\zeta_t)_{t\geq 0}\quad{for\,\,\, each} \quad c>0.$$

\noindent
(3) Second moments of increments: for $0\leq s <t$, 

$$E(\zeta_t-\zeta_s)^2\leq C|t-s|^2,$$
if $s,t\leq T$, for any $T>0$ with $C=C(T)$;
$$E(\zeta_t-\zeta_s)^2\geq C|t-s|^2,$$
if $s,t\geq\varepsilon$, for any $\varepsilon>0$ with $C=C(\varepsilon)$;
$$E(\zeta_t-\zeta_s)^2\geq C|t-s|^h,$$
for all $s,t\geq 0$.

\noindent
(4) Positive correlation of increments: For $0\leq r <v$, $0\leq s <t$,
$$
E((\zeta_v-\zeta_r)(\zeta_t-\zeta_s))>0.
$$

\noindent
(5) Long-range dependence: For $0\leq r<v\leq s<t$,
$$\lim_{T\to\infty}T^{3-h}E((\zeta_v-\zeta_r)(\zeta_{t+T}-\zeta_{s+T}))=\frac{h(h-1)(h-2)^2}{2}(t-s)(v^2-r^2).$$%

\noindent
(6) $\zeta$ is not Markov.
\end{proposition}

We omit the proof of this proposition, as it is easily obtained from 
(\ref{eq:3.1}) and (\ref{eq:3.2}). Integral representations of sfBm are given in \cite{BGT2}. Analogous results can be obtained for nsfBm from (\ref{eq:3.2}) and the known integral representations of fBm (see Remark 3.3 (c)).

Let us recall the function $R$ defined by (\ref{eq:1.7}):
\begin{equation}
\label{eq:3.5p}
R(s,t)=-\biggl(s^2 \log s+t^2 \log t-\frac{1}{2}[(s+t)^2 \log (s+t) +(s-t)^2 \log |s-t|]\biggr), \quad s, t\geq 0;
\end{equation}
we already know that it is positive-definite.

The proximity of the process $\eta =(\eta_t)_{t\geq 0}$ with covariance (\ref{eq:3.5p}) to nsfBm is seen if we look at formula (\ref{eq:3.4}), since it easy to verify that (\ref{eq:3.5p}) can be written as 
\begin{equation}
\label{eq:3.5}
R(s,t)=2\int^{s\wedge t}_0\int^s_r\int^t_r(u+u'-2r)^{-1}du'dudr.
\end{equation}

It is also worthwhile to recall that nsfBm is related to a particle system for $d<\alpha, \eta$ corresponds to the borderline case $d=\alpha$, and the case $\alpha<d<2\alpha$ leads to sfBm (see the Introduction).

Some properties of $\eta$, such as self-similarity and long-range dependence, are analogous to those of nsfBm  and are easy to obtain. Path continuity follows for example from \cite{BGT6}. 
The long-range dependence is given by
$$
\lim_{T\rightarrow \infty} TE((\eta_v-\eta_r)(\eta_{t+T}-\eta_{s+T}))=(t-s)(v^2-r^2),
$$
$0\leq r<v\leq s <t$.  The most interesting property of $\eta$ is the fact that, unlike nsfBm, $\eta$ it not a semimartingale, so, in this sense $\eta$ is closer to sfBm. We state this property as the next proposition.

\begin{proposition}
\label{P:3.5}
The process $\eta $ is not a semimartingale.
\end{proposition}

\subsection{Proofs}

\vglue .25cm
\noindent
{\bf Proof of Theorem 3.1} Existence of the process $\vartheta$ follows from Remark 3.3 (c).

From (\ref{eq:3.1}) we have 
$$K(s,t)=\frac{1}{2}h(h-1)(h-2)\int^t_0\int^s_0 K_0(u,v)dudv,$$
where $K_0$ is given by (\ref{eq:3.3}). Hence the theorem is proved. \hfill $\Box$

\vglue.5cm
It is easy to show that $K_0$ is not positive-definite for $h>4$, so the bfBm with $H=\frac{1}{2}$ cannot be extended to $k>2$.
\vglue.5cm
\noindent
{\bf Proof of Proposition 3.5} We will show that $\eta$ has infinite variation on $[0,1]$, and its quadratic variation is zero.

From (\ref{eq:3.5}) for $s<t$, we obtain by elementary calculus 
\begin{equation}
\label{eq:3.7}
E(\eta_t-\eta_s)^2=\int^{t-s}_0\int^u_0\left[\log(2s+u+u')-\log(u-u')\right]du'du.
\end{equation}
We put $s=k/n,t=(k+1)/n, e^4\leq k\leq n$. Using 
$\log (2k/n)\leq \log (2k/n+u+u')\leq \log ((2k+2)/n)$ for $0\leq u,u'\leq 1/n$, it is not hard to see that (\ref{eq:3.7}) implies
\begin{equation}
\label{eq:3.8}
\frac{1}{2n^2}\log k\leq E(\eta_{(k+1)/n}  -\eta _{k/n})^2\leq\frac{1}{n^2}\log k.
\end{equation}
(The assumption $k\geq e^4$ is used in the upper estimate, where we need to have that $(\log(2k+2)+2)/2\leq \log k$. The lower estimate holds for any $k\geq 1$.)

From the upper estimate in (\ref{eq:3.8}) we infer that
$$\lim_{n\to\infty}\sum^{n-1}_{k=0}E(\eta_{(k+1)/n}-\eta_{k/n})^2=0,$$
hence the quadratic variation of $\eta$ cannot be different from $0$.

By the Gaussian property of $\eta$, (\ref{eq:3.8}) implies that
\begin{equation}
\label{eq:3.9}
E|\eta_{(k+1)/n}-\eta_{k/n}|\geq C\frac{(\log k)^{1/2}}{n}.
\end{equation}

To show that the variation of $\eta$ is infinite, we repeat the argument of the proof of Lemma 2.1 in \cite{BGT1}. By the $0-1$ law for Gaussian processes, the variation on $[0,1]$ is either finite a.s. or infinite a.s.  If it were finite a.s. (we denote it by $|\eta|_1$), then by the Fernique theorem \cite{F} we would have $E|\eta|_1<\infty$. On the other hand (\ref{eq:3.9}) implies that
$$E|\eta|_1\geq \frac{1}{n}\sum^{n-1}_{k=1}
(\log k)^{1/2}\to\infty\quad{\rm  as}\quad n\to\infty .$$
Therefore $|\eta|_1=\infty$ a.s. \hfill $\Box$

\vglue.25cm
\noindent
{\bf Acknowledgment} We thank the hospitality of the Institute of Mathematics, National University of Mexico (UNAM), where this paper was written.

\end{document}